\documentclass{ifacconf}

\usepackage{graphicx}     
\usepackage{natbib} 

\usepackage{amsmath}
\usepackage{amsfonts}
\usepackage{units}
\usepackage{csquotes}
\usepackage{color}
\usepackage{subcaption}
\usepackage{mathptmx}
\allowdisplaybreaks
\begin{document}
\begin{frontmatter}

\title{Suppression of Oscillations in Two-Class Traffic by Full-State Feedback \thanksref{footnoteinfo}} 

\thanks[footnoteinfo]{Mark Burkhardt acknowledges financial support through the Baden-W\"urttemberg Stipendium of the Baden-W\"urttemberg Stiftung. }

\author[First]{Mark Burkhardt} 
\author[Second]{Huan Yu} 
\author[Second]{Miroslav Krstic}

\address[First]{Institute for System Dynamics, University of Stuttgart, Stuttgart, Germany (e-mail: mark.burkhardt@isys.uni-stuttgart.de)}
\address[Second]{Department of Mechanical and Aerospace Engineering, University of California, San Diego, USA (e-mail: \{huy015,krstic\}@ucsd.edu)}

\begin{abstract}          
This paper develops a full-state feedback controller that damps out oscillations in traffic density and traffic velocity whose dynamical behavior is governed by the linearized two-class Aw-Rascle (AR) model. Thereby, the traffic is considered to be in the congested regime and subdivided in two classes whereas each class represents vehicles with the same size and driver's behavior. The macroscopic second-order two-class AR model consists of four first order hyperbolic partial differential equations (PDEs) and introduces a concept of area occupancy to depict the mixed density of two-class vehicles in the traffic. Moreover, the linearized model equations show heterodirectional behavior with both positive and negative characteristic speeds in the congested regime. The control objective is to achieve convergence to a constant equilibrium in finite time. The control input is realized by ramp metering acting at the outlet of the considered track section. The backstepping method is employed to design full-state feedback for the $4\times 4$ hyperbolic PDEs. The performance of the full-state feedback controller is verified by simulation.  
\end{abstract}
\vspace{-7pt}
\end{frontmatter}
%===============================================================================

\section{Introduction}
\vspace{-8pt}
A common issue in everyday life is congested traffic. Since more and more people own a car, highways become more crowded leading to a problem of growing impact. A well-known phenomenon that occurs in congested traffic is stop and go traffic. In fact, stop and go traffic is characterized by oscillations in traffic density and velocity causing higher fuel consumption and a higher risk of accidents. \\
Traffic models can be employed to investigate stop and go traffic as well as developing solutions to avoid this phenomenon. They are categorized in microscopic, mesoscopic and macroscopic models. Since macroscopic models describe traffic as a distributed parameter system in traffic densities, flows or velocities, they are more suitable to investigate density and velocity oscillations. Their model equations are governed by PDEs in traffic densities and velocities along an investigated track section. Moreover, macroscopic traffic models either capture homogeneous traffic or heterogeneous traffic. Homogeneous traffic consists of vehicles whose size and driver's behavior are the same, whereas heterogeneous considers traffic containing different classes which are defined as a group of vehicles with the same properties, see~\cite{logghe2003dynamic}. For instance, motorcycles and trucks as well as tolerant and aggressive drivers can be distinguished in heterogeneous traffic yielding a different dynamical behavior. Notice that a macroscopic traffic model that covers heterogeneous traffic and hence introduces multiple classes of vehicles is denoted as macroscopic multi-class traffic model. First-order macroscopic multi-class traffic models, like~\cite{wong2002multi},~\cite{van2008fastlane},~\cite{fan2015heterogeneous}, and second-order models, e.g.~\cite{SOMC_Gupta},~\cite{SOMC_ExtendedSG},~\cite{SOMC_TangDerviedfromCarfollow}, are distinguished in literature depending on whether a second PDE is introduced to model the velocity dynamics of each class. The main focus of this paper is on the extended AR model presented in~\cite{MCAR} which is a nonlinear second-order multi-class traffic model. Therein, the concept of area occupancy is introduced yielding a coupling between the vehicle classes. In this work, traffic contains two classes. Thus, the extended AR model is evaluated for two classes and introduced as the two-class AR traffic model. \\
Control of macroscopic traffic models is addressed in previous work. For instance,~\cite{karafyllis2018feedback} introduces a macroscopic second-order traffic model consisting of two PDEs and develops a stabilizing boundary feedback law only depending on the inlet speed. The corresponding control law is deduced by formulating a boundary condition for the characteristic form of the model. Furthermore,~\cite{bekiaris2019feedback} considers an Aw-Rascle-Zhang-type model describing traffic that is equipped with Adaptive Cruise Control yielding a spatially distributed input. The corresponding linearized system is unstable, but is stabilized by a feedback law that eliminates the source term. In addition, control design for multi-class traffic is carried out in~\cite{deo2009model},~\cite{liu2016model} and~\cite{pasquale2015two}. Whereas~\cite{deo2009model} focuses on a model predictive control approach for a macroscopic traffic flow model extended to multi-class flows,~\cite{pasquale2015two} presents an optimal control problem and numerical solution algorithm with respect to the total emissions and total time spent for a traffic model distinguishing between two classes. Moreover,~\cite{liu2016model} develops a model predictive controller for multi-class traffic and emission models. \\
Traffic management systems like ramp metering or variable speed limits can be employed to damp the density and velocity oscillations yielding evenly distributed traffic on freeways, see e.g.~\cite{VSLARZ},~\cite{ARZYU},~\cite{ZhangPrieur_ExpStability_Positivehyp}. The overall goal of this work is to develop a full state feedback controller in order to damp out traffic density and velocity oscillations in finite time. Similar to a valve at the end of a pipe, a ramp metering system is thus assumed to be installed at the end of the investigated track section yielding a boundary control problem. In this work, the backstepping technique, see e.g.~\cite{Deutscher_2x2_2},~\cite{DiMeglio_Control}, is employed to achieve the damping in finite time. The control design is carried out for a heterodirectional system of $3+1$ transport PDEs corresponding to a special case of the results presented in~\cite{HuKrsticGeneralPDEs}. Thereby, $3+1$ transport PDEs depicts three PDEs with downstream convection and one PDE with upstream convection. The transport velocities of the transport systems are addressed by the term characteristic speeds in the following.  \\
Contribution of this paper: this work yields the first result on boundary feedback control of a macroscopic second-order multi-class traffic model using backstepping. Thus, it contributes to the application of a theoretical control design method to realistic traffic model. Moreover, it contributes to traffic modeling by investigating the characteristic speeds of the linearized two-class AR model and identifying that only one characteristic speed is negative in its congested regime. \\
The paper is structured as follows: Section $2$ introduces the model equations of the two-class AR model. Furthermore, it covers the linearization around a constant equilibrium and the analysis of the characteristic speeds for the resulting linearized model equations. In a last step, this section presents a transformation into Riemann coordinates. Afterwards, Section $3$ deals with the full state feedback design. Therein, the control objective is formulated, the desired control law is developed and the result is summarized in a theorem. Simulation results in order to verify the controller performance are given in Section $4$, followed by concluding remarks in Section $5$.
\section{Problem statement}
\vspace{-5pt}
This section presents the two class AR model that models traffic consisting of two vehicle classes and the control design model. First, the nonlinear model equations are introduced. Afterwards, the equations are linearized around a constant equilibrium yielding the linearized two class AR model. In the next step, the boundary conditions are discussed and linearized. Finally, the signs of the characteristic speeds are investigated and the control design model is deduced by introducing a transformation into Riemann coordinates.
\vspace{-1pt} 
\subsection{Two-class AR model}
\vspace{-3pt}
The Extended AR model in~\cite{MCAR} is formulated for an arbitrary amount of vehicle classes. In this work, traffic distinguishing two vehicle classes is considered yielding the two class AR model
\begin{subequations}
	\begin{align}
	\partial_t \rho_1 \hspace{-2pt} =&\hspace{-1pt} -\partial_x(\rho_1 v_1), \label{sec2:eq:TwoClassAR_nonlin1} \\
	\partial_t (v_1\hspace{-2pt}+\hspace{-2pt}p_1(AO))\hspace{-2pt}+\hspace{-2pt}v_1\partial_x (v_1\hspace{-2pt}+\hspace{-2pt}p_1(AO))\hspace{-2pt}=&\frac{\hspace{-1pt}V_{e,1}(AO)-\hspace{-1pt}v_1}{\tau_1} , \label{sec2:eq:TwoClassAR_nonlin2} \\
	\partial_t \rho_2  \hspace{-2pt}=&\hspace{-1pt} -\partial_x(\rho_2 v_2), \label{sec2:eq:TwoClassAR_nonlin3}\\
	\partial_t (v_2\hspace{-2pt}+\hspace{-2pt}p_2(AO))\hspace{-2pt} +\hspace{-2pt}v_2\partial_x (v_2+p_2(AO))\hspace{-2pt}=&\frac{\hspace{-1pt}V_{e,2}(AO)-\hspace{-1pt}v_2}{\tau_2}. \label{sec2:eq:TwoClassAR_nonlin4}
	\end{align}
	\label{sec2:eq:TwoClassAR}
	\vspace{-5pt}
\end{subequations} \\
The initial conditions are defined as $\rho_i(x,0) = \rho_{i,0}(x) \in \mathcal{L}^\infty([0,L])$ as well as $v_i(x,0) = v_{i,0}(x) \in \mathcal{L}^\infty([0,L])$ and the boundary conditions~\eqref{sec2:eq:Nonlinear_BCs} are discussed later on. The traffic density and velocity of vehicle class $i$ are denoted by $\rho_i(x,t)$ and $v_i(x,t)$ on the domain $(x,t) \in (0,L)\times(0,\infty)$, whereas $L$ is the length of the considered track section. In general, the traffic density is the amount of vehicles per unit length. The four nonlinear first-order hyperbolic PDEs~\eqref{sec2:eq:TwoClassAR} are coupled and represent conservation laws. The terms on the right hand side of~\eqref{sec2:eq:TwoClassAR_nonlin2} and~\eqref{sec2:eq:TwoClassAR_nonlin4} describe the adaption of the vehicles to their desired velocity in adaption time $\tau_i$. According to the simplified definition in~\cite{MCAR}, the area occupancy measurement is defined as
\begin{equation}
\label{sec2:eq:AO_def}
AO(\rho_1,\rho_2)=\frac{a_1L\rho_1+a_2L\rho_2}{WL}
\end{equation}
and thus depends on both traffic densities. The parameter $a_i$ is the surface that is occupied by one vehicle of a class. Additionally, $W$ is the width of the investigated track section. Physically interpreted, $AO$ is the percentage of road space that is occupied if vehicles of class one are distributed with density $\rho_1(x,t)$ and vehicles of class two are distributed with density $\rho_2(x,t)$ along the considered track section. The traffic pressure
\begin{equation}
\label{sec2:eq:Pressure_AO}
p_i(AO) = V_{i}\left(\frac{AO(\rho_1,\rho_2)}{\overline{AO}_{i}}\right)^{\gamma_i}
\end{equation}
and the equilibrium speed $AO$ relationship
\begin{equation}
\label{sec2:eq:Equil_speed_AO}
V_{e,i}(AO) = V_{i}\left(1-\left(\frac{AO(\rho_1,\rho_2)}{\overline{AO}_{i}}\right)^{\gamma_i}\right)
\end{equation}
depend on the area occupancy. The traffic pressure models the reaction to other vehicles. Qualitatively, it holds that a more crowded freeway yields higher $AO$ leading to a higher experienced traffic pressure. The equilibrium speed $AO$ relationship describes the desired velocities of the vehicles and is based on the model of Greenshield~\cite{Greenshield}. Moreover, both functions introduce three additional parameters: the free-flow velocity $V_i>0$, traffic pressure exponent $\gamma_i>1$ and maximum area occupancy $0<\overline{AO}<1$. The free-flow velocity $V_i$ represents the desired velocity without taking other vehicles into account, the traffic pressure exponent $\gamma_i$ is a degree of freedom to adjust the pressure function to realistic data and the maximum area occupancy $\overline{AO}_i$ is the area occupancy value that implies whether a vehicle class is jammed, i.e. desired velocity $V_{e,i}(\overline{AO}_i) = 0$. In fact, it holds that $V_{e,i}(AO) = V_i-p_i(AO)$. Thus, a more crowded freeway implying higher traffic pressure yields a more reduced desired velocity. Besides, notice that $V_{e,i}(0)=V_i$ which equals to the definition of $V_i$.
\vspace{-9pt}
\subsection{Linearized two-class AR model}
\vspace{-5pt}
Next, the two class AR model is linearized around the constant steady state $z^*=(\rho_1^*,v_1^*,\rho_2^*,v_2^*)^T$. Notice that the steady state satisfies
\begin{equation}
v_i^*(\rho_1^*,\rho_2^*) = V_{e,i}(AO(\rho_1^*,\rho_2^*)) \label{sec2:eq:EquilibriumState_cond}
\end{equation}
which is obtained after inserting $z^*$ in~\eqref{sec2:eq:TwoClassAR}. According to~\eqref{sec2:eq:EquilibriumState_cond}, the steady state velocities follow, if the steady state densities are determined. The perturbations around the steady state are defined as
\begin{equation}
\tilde{\rho}_i(x,t) = \rho_i(x,t)-\rho_i^*, \quad \tilde{v}_i(x,t)= v_i(x,t)-v_i^*  \label{sec2:eq:PerturbationsDefinition}
\end{equation}
and are summarized in the state vector $z=(\tilde{\rho}_1,\tilde{v}_1,\tilde{\rho}_2,\tilde{v}_2)^T$. Consequently, the linearized two-class AR model is given by
\begin{equation}
\label{sec2:eq:LinearizedModelEq}
J_tz_t+J_xz_x + Jz = 0,
\end{equation}
with Jacobian matrices  
\begin{align}
J_t &= \left[\begin{array}{cccc}
1 & 0 & 0 & 0 \\
\beta_{11} & 1 & \beta_{12} & 0 \\
0 & 0 & 1 & 0 \\
\beta_{21} & 0 & \beta_{22} & 1 \\
\end{array}\right], \text{ } J_x = \left[\begin{array}{cccc}
v_1^* & \rho_1^* & 0 & 0 \\
v_1^*\beta_{11} & v_1^* & v_1^* \beta_{12} & 0 \\
0 & 0 & v_2^* & \rho_2^* \\
v_2^*\beta_{21} & 0 & v_2^*\beta_{22} & v_2^* \\
\end{array}\right], \notag \\
J&=\left[\begin{array}{cccc}
0 & 0 & 0 & 0 \\
\frac{1}{\tau_1}\beta_{11} & \frac{1}{\tau_1} & \frac{1}{\tau_1} \beta_{12} & 0 \\
0 & 0 & 0 & 0 \\
\frac{1}{\tau_2}\beta_{21} & 0 & \frac{1}{\tau_2}\beta_{22}   & \frac{1}{\tau_2} \\
\end{array}\right]
\end{align}
including the abbreviation
\begin{equation}
\beta_{ij}(\rho_1^*,\rho_2^*) = \left.\frac{\partial p_i(AO(\rho_1,\rho_2))}{\partial \rho_j}\right|_{\rho_1 = \rho_1^*,\rho_2 = \rho_2^*},
\end{equation}
where $i,j=1,2$. \\
Four boundary conditions are discussed in order to complete the linearized two-class AR model. First, it is assumed that the incoming traffic flow has fixed densities. Moreover, the remaining two conditions describe that the overall traffic flow entering and leaving the investigated track section is constant. Thus, the boundary conditions are given by
\begin{subequations}
\begin{align}
	&\rho_i(0,t) = \rho_i^*, q_1(0,t)+q_2(0,t)=\rho_1^*v_1^*+\rho_2^*v_2^*,\\
	& q_1(L,t)+q_2(L,t) =\rho_1^*v_1^*+\rho_2^*v_2^*+U(t).
\end{align}
\vspace{-7pt}
\label{sec2:eq:Nonlinear_BCs}
\end{subequations} \\
Since the traffic flow is defined as $q_i=\rho_iv_i$, the latter conditions lead to nonlinear equations in $\rho_i$ and $v_i$. The control input $U(t)$ is realized by a ramp metering at the outlet of the track section regulating the leaving traffic flow similar to a valve at the end of a pipe. Notice that $U(t)$ describes the perturbation of the traffic flow caused by the ramp. The linearized boundary conditions then become
\begin{subequations}
	\begin{align}
	0 & \hspace{-2pt}=\hspace{-2pt} \tilde{\rho}_i(0,t), \label{sec2:eq:BC0_lin_Dens1} \\
	0 & \hspace{-2pt}=\hspace{-2pt} v_1^*\tilde{\rho}_1(0,t)\hspace{-2pt}+\hspace{-2pt}\rho_1^*\tilde{v}_1(0,t)\hspace{-2pt}+\hspace{-2pt}v_2^*\tilde{\rho}_2(0,t)\hspace{-2pt}+\hspace{-2pt}\rho_2^*\tilde{v}_2(0,t),  \label{sec2:eq:BC0_lin_Flow}\\
	U(t) \hspace{-2pt}&= \hspace{-2pt} v_1^*\tilde{\rho}_1(L,t)\hspace{-2pt}+\hspace{-2pt}\rho_1^*\tilde{v}_1(L,t)\hspace{-2pt}+\hspace{-2pt}v_2^*\tilde{\rho}_2(L,t)\hspace{-2pt}+\hspace{-2pt}\rho_2^*\tilde{v}_2(L,t). \label{sec2:eq:BCL_lin_Flow_wContrl} 
	\end{align}
	\label{sec2:eq:BC_lin}
	\vspace{-14pt}
\end{subequations}
\subsection{Free-flow/congested regime analysis}
Based on the steady state around which the linearization is carried out and the parameters occurring in the PDEs, the linearized two class AR model captures two fundamentally different dynamical behaviors: either the described traffic is in the free-flow regime or in the congested regime. In the free-flow regime, all characteristic speeds are positive and there is only information propagating downstream. This behavior is clarified as homo-directional behavior because the information travels in one direction. However, traffic in the congested regime is characterized by information propagating upstream, i.e. heterodirectional behavior. In the following, the amount of negative characteristic speeds in the latter regime is investigated. Notice that stop and go traffic only occurs in the congested regime, thus analyzing the dynamical behavior in this regime is crucial for the control design. The characteristic speeds
\begin{equation}
\lambda_i = v_i^*, \text{ } i=1,2, \quad \lambda_{3/4} = \frac{v_1^*+v_2^*-\beta_{11}\rho_1^*-\beta_{22}\rho_2^*\pm\Delta}{2}, \label{sec2:eq:CharSpeed2}
\end{equation}
where
\begin{equation}
\Delta(\rho_1^*,\rho_2^*) \hspace{-2pt}=\hspace{-3pt} \sqrt{\hspace{-2pt}\left(\hspace{-1pt}\beta_{22}\rho_2^*\hspace{-2pt}-\hspace{-1pt}\beta_{11}\rho_1^*\hspace{-2pt}+\hspace{-1pt}v_1^*\hspace{-2pt}-\hspace{-1pt}v_2^*\right)^2\hspace{-4pt}+\hspace{-1pt}4\beta_{11}\beta_{22}\rho_1^*\rho_2^*}\hspace{-1pt}
\end{equation}
are computed by calculating the eigenvalues of the Jacobian $J_t^{-1}J_x$. In fact,~\cite{AnalyseARZ_CharSpeed} shows that the relation 
\begin{equation}
\label{sec2:eq:CharSpeed_Order}
\lambda_4 \leq \min\{\lambda_1,\lambda_2\} \leq \lambda_3 \leq \max\{\lambda_1,\lambda_2\} 
\end{equation}
holds. The characteristic speeds $\lambda_1>0$ and $\lambda_2>0$ are the steady state velocities which are positive due to model validity. Furthermore, $\lambda_{3}$ correspond to the traffic flow that is caused by the vehicles overtaking each other yielding a positive sign. Thus, $\lambda_4$ is the only characteristic speed that may be smaller than zero. In the following, the traffic captured by the linearized two-class AR model~\eqref{sec2:eq:LinearizedModelEq} is in the congested regime, if
\begin{equation}
\label{sec2:eq:CongestedRegime_Def}
	\lambda_1,\lambda_2,\lambda_3>0, \text{ } \lambda_4<0
\end{equation}
and the traffic is in free-flow regime if
\begin{equation}
	\lambda_1,\lambda_2,\lambda_3,\lambda_4>0.
\end{equation}
Throughout the rest of this paper, it is assumed that the steady state densities and parameters are chosen such that~\eqref{sec2:eq:CongestedRegime_Def} is satisfied.
\vspace{-3pt}
\subsection{Control design model}
\vspace{-3pt}
Next, the linearized two-class model~\eqref{sec2:eq:LinearizedModelEq} with~\eqref{sec2:eq:BC_lin} is transformed into Riemann coordinates yielding the control design model. The preparation covers two goals. First, the characteristic form of the PDEs is deduced. Second, the states are sorted such that the characteristic speeds occur in ascending order and the diagonal elements of the source term become zero in order to keep the computations concise and increase the readability by adapting to the notation of~\cite{HuKrsticGeneralPDEs}. Notice that an unique ascending order of the characteristic speeds is defined as soon as one of the steady state velocities is larger than the other. In this work, $v_1^*>v_2^*$ is assumed and thus the first vehicle class corresponds to faster vehicles. The anticipated order of characteristic speeds is then given by $\lambda_4<0<\lambda_2<\lambda_3<\lambda_1$. \\
For that reason, the transformation 
\begin{equation}
\label{sec3:eq:SummarizedTransformation}
w_c = \left[\begin{array}{cccc}
0 & e^{-\frac{\hat{J}_{22}}{v_2^*}x} & 0 & 0 \\
0 & 0 & e^{-\frac{\hat{J}_{33}}{\lambda_3}x} & 0 \\
e^{-\frac{\hat{J}_{11}}{v_1^*}x} & 0 & 0 & 0 \\
0 & 0 & 0 & e^{-\frac{\hat{J}_{44}}{\lambda_4}x}
\end{array}\right]
\Theta^{-1}z
\end{equation}
with $\Theta$ such that
\begin{equation}
\text{diag}(\lambda_1,\lambda_2,\lambda_3,\lambda_4) = \Theta^{-1}J^{-1}_tJ_x\Theta
\end{equation}
and $\Theta = \left\{\theta_{ij}\right\}_{1\leq i\leq 4,1\leq j \leq 4}$ as well as $\hat{J}=-\Theta^{-1}J^{-1}_tJ\Theta$ with $\{\hat{J}_{ij}\}_{1 \leq i \leq 4, 1\leq j \leq 4}$ is introduced. The Riemann coordinates are denoted as $w_c=(w_1,w_2,w_3,w_4)^T$ and in order to keep the computations concise, the first three coordinates are summarized in $w=(w_1,w_2,w_3)^T$. This transformation is applied to the system with decoupled partial time derivatives, i.e. the model equations that result from multiplying~\eqref{sec2:eq:LinearizedModelEq} with $J_t^{-1}$. The transformed model equations are given by
\begin{subequations}
	\begin{align}
	w_t+\Lambda^+
	w_x&= \Sigma^{++}(x)
	w + \Sigma^{+-}(x)w_4, \label{sec3:eq:ControlDesignModel_w123PDE} \\
	w_{4t}- \Lambda^- w_{4x} &= \Sigma^{-+}(x) w \label{sec3:eq:ControlDesignModel_w4PDE}
	\end{align}
	\label{sec3:eq:ControlDesignModel}
\end{subequations} \\
with
\begin{align}
\Lambda^+ &= \text{diag}(v_2^*,\lambda_3,v_1^*),\text{ } \Lambda^-= -\lambda_4, \\
\Sigma^{++}(x) &= \left[\begin{array}{ccc}
0 & \bar{J}_{12}(x) & \bar{J}_{13}(x) \\
\bar{J}_{21}(x) & 0 & \bar{J}_{23}(x) \\
\bar{J}_{31}(x)& \bar{J}_{32}(x) & 0 \\ 
\end{array}\right], \\
\Sigma^{+-}(x) &= \left[\begin{array}{ccc}
\bar{J}_{14}(x) & \bar{J}_{24}(x) & \bar{J}_{34}(x)
\end{array}\right]^T,\\
\Sigma^{-+}(x) &= \left[\begin{array}{ccc}
\bar{J}_{41}(x) & \bar{J}_{42}(x) & \bar{J}_{43}(x)
\end{array}\right]. 
\end{align}
Since the traffic is in congested regime and thus $\lambda_4<0$,~\eqref{sec3:eq:ControlDesignModel_w123PDE} captures the information propagating downstream and~\eqref{sec3:eq:ControlDesignModel_w4PDE} the information propagating upstream. The coefficients $\bar{J}_{ij}$ are omitted here. They are bounded and do not change their sign on the whole domain of $x$. As a remark, the input of the transformed system is given by
\begin{equation}
\label{sec3:eq:SummarizedTransformation_Ubar}
\bar{U}(t) = e^{-\frac{\hat{J}_{44}}{\lambda_4}L}\frac{1}{\kappa}U(t)
\end{equation}
with the abbreviation $\kappa = v_1^*\theta_{14}+\rho_1^*\theta_{24}+v_2^*\theta_{34}+\rho_2^*\theta_{44}$. The numerical investigations that were considered while carrying out this work, show $\kappa \neq 0$ if $v_1^*,v_2^*>0$.  \\
The transformation~\eqref{sec3:eq:SummarizedTransformation} also needs to be applied to the boundary conditions~\eqref{sec2:eq:BC_lin}. For that reason, the boundary conditions are formulated with the state vector $z$ and afterwards the transformation law is inserted yielding the conditions
\begin{subequations}
	\begin{align}
	&w(0,t) = \bar{Q}_0w_4(0,t), \label{sec3:eq:ControlDesignModel_bc0} \\
	&w_4(L,t) = \bar{R}_1w(L,t)+\bar{U}(t). \label{sec3:eq:ControlDesignModel_bcL}
	\end{align}
	\label{sec3:eq:ControlDesignModel_bc}
	\vspace{-9pt}
\end{subequations} \\
The matrices $\bar{Q}_0$ and $\bar{R}_1$ are straightforward to compute and thus omitted due to space constraints. \\
In the following,~\eqref{sec3:eq:ControlDesignModel} and~\eqref{sec3:eq:ControlDesignModel_bc} are considered as control design model. The performed transformation is invertible and therefore the linearized two-class AR model and the control design model have the same stability properties.
\section{Full state feedback control design}
\vspace{-5pt}
In the following, a full-state feedback control design for the system of four coupled hyperbolic PDEs given by~\eqref{sec3:eq:ControlDesignModel} with boundary conditions~\eqref{sec3:eq:ControlDesignModel_bc} is carried out in order to achieve finite time convergence to zero for initial conditions $w_j(x,0) \in \mathcal{L}^\infty[0,L]$. The overall goal is to damp out stop-and-go traffic in the congested regime and achieve convergence to the constant steady state in a finite time. The term of stop-and-go traffic refers to oscillations of the density and velocity perturbations around their constant equilibrium values along the highway. The full-state feedback controller is designed by applying the backstepping control design, see~\cite{HuKrsticGeneralPDEs}. The states of the target system are denoted as $(\alpha,\beta)^T$, where $\alpha=(\alpha_1,\alpha_2,\alpha_3)^T$. The kernels of the backstepping transformation are denoted by $K(x,\xi)=\left\{k_{1j}(x,\xi)\right\}_{1\leq j \leq 3}$ and $L_{11}(x,\xi)$. Then, the backstepping transformation is defined as
\begin{subequations}
	\begin{align}
	\alpha_j(x,t)\hspace{-2pt}&=\hspace{-2pt} w_j(x,t),\quad j=1,2,3, \label{sec4:eq:Cntrl_BacksteppingTrafo1} \\
	\beta(x,t)\hspace{-2pt} &=\hspace{-2pt} w_4(x,t)\hspace{-2pt}-\hspace{-4pt}\int_0^x\hspace{-4pt}\left(\hspace{-1pt}K(x,\xi)\hspace{-1pt} w(\xi,t)\hspace{-1pt}+\hspace{-1pt} L_{11}(x,\xi)w_4(\xi,t)\hspace{-1pt}\right)\hspace{-1pt}d\xi \label{sec4:eq:Cntrl_BacksteppingTrafo2}.
	\end{align}
	\label{sec4:eq:Cntrl_BacksteppingTrafo}
	\vspace{-7pt}
\end{subequations} \\
Notice that $K(x,\xi)$ and $L_{11}(x,\xi)$ are defined on a triangular domain $\mathcal{T} = \{0\leq \xi\leq x\leq 1\}$.
Furthermore, the choice of the well-posed target system is 
\begin{subequations}
	\begin{align}
	&\alpha_t+\Lambda^+\alpha_x =\Sigma^{++}(x)
	\alpha+ \Sigma^{+-}(x)\beta \notag \\ 
	&+ \int_0^x C^+(x,\xi)\alpha(\xi,t)d\xi+ \int_0^x C^-(x,\xi)\beta(\xi,t) d\xi, \label{sec4:eq:Cntrl_TargetSystem_a123PDE} \\
	&\beta_t= \Lambda^-\beta_x \label{sec4:eq:Cntrl_TargetSystem_bPDE}.
	\end{align}
	\label{sec4:eq:Cntrl_TargetSystem}
	\vspace{-11pt}
\end{subequations} \\
The coefficients $C^+(x,\xi) \in \mathbb{R}^{3\times 3}$ and $C^-(x,\xi)\in \mathbb{R}^{3 \times 1}$ are defined on the same triangular domain $\mathcal{T}$ and are determined later on. Besides, the boundary conditions of the target system are
\begin{subequations}
	\begin{align}
	\alpha(0,t) &= \bar{Q}_0 \beta(0,t), \label{sec4:eq:Cntrl_TargetSystem_BC0} \\
	\beta(L,t) &= 0. \label{sec4:eq:Cntrl_TargetSystem_BCL}
	\end{align}
	\label{sec4:eq:Cntrl_TargetSystem_BC}
	\vspace{-11pt}
\end{subequations} \\
The target system~\eqref{sec4:eq:Cntrl_TargetSystem} with~\eqref{sec4:eq:Cntrl_TargetSystem_BC} converges to its equilibrium  at zero
\begin{equation}
\alpha_{e,j}(x)\equiv \beta_e(x) \equiv 0, \quad j=1,2,3, \text{ } t\geq 0, \text{ } x\in [0,L]
\end{equation}
in the finite time
\begin{equation}
\label{sec4:eq:FiniteConvergenceTime_tF}
t_F = \frac{L}{v_2^*}+\frac{L}{-\lambda_4}.
\end{equation}
The proof is given in Lemma $3.1$ in~\cite{HuKrsticGeneralPDEs}. Differentiating~\eqref{sec4:eq:Cntrl_BacksteppingTrafo2} with respect to space and time, inserting the resulting derivatives and~\eqref{sec3:eq:ControlDesignModel_bc0} in~\eqref{sec4:eq:Cntrl_TargetSystem_bPDE} yields kernel equations that determine $K(x,\xi)$ and $L_{11}(x,\xi)$. Employing the method of characteristics afterwards yields
\begin{align}
L_{11}(x,\xi) =& -\frac{1}{\lambda_4}K(x-\xi,0)\Lambda^+\bar{Q}_0 \notag \\
+\int_0^{-\frac{\xi}{\lambda_4}}&K(\lambda_4\nu+x,\lambda_4\nu+\xi)\Sigma^{+-}(\lambda_4\nu+\xi)d\nu \label{sec4:eq:kerneleq_PDEL_Linserted}
\end{align}
and three coupled first order hyperbolic PDEs with three boundary conditions
\begin{subequations}
	\begin{align}
	0=&\lambda_4K_x(x,\xi)+\Lambda^+K_\xi(x,\xi)+K(x,\xi)\Sigma^{++}(\xi) \notag \\
	&-\frac{1}{\lambda_4}K(x-\xi,0)\Lambda^+\bar{Q}_0\Sigma^{-+}(\xi) \notag \\
	+&\hspace{-2pt}\int_0^{-\frac{\xi}{\lambda_4}}\hspace{-2pt}K(\lambda_4\nu\hspace{-2pt}+\hspace{-2pt}x,\lambda_4 \nu\hspace{-3pt}+\hspace{-2pt}\xi)\Sigma^{+-}(\lambda_4\nu\hspace{-2pt}+\hspace{-2pt}\xi)d\nu\Sigma^{-+}(\xi) \label{sec4:eq:kerneleq_PDEK_Linserted} \\
	0=&K(x,x)\Lambda^++\Lambda^-K(x,x)+\Sigma^{-+}(x). \label{sec4:eq:kerneleq_BCxx_Linserted}
	\end{align}
	\label{sec4:eq:kerneleq}
	\vspace{-11pt}
\end{subequations} \\
As shown in Theorem $3.3$ of~\cite{HuKrsticGeneralPDEs}, the kernel equations~\eqref{sec4:eq:kerneleq_PDEL_Linserted} and~\eqref{sec4:eq:kerneleq} are a well-posed system of equations and thus there exist unique solutions $K(x,\xi)$ and $L_{11}(x,\xi)$ in $L\infty(\mathcal{T})$. Furthermore, differentiating~\eqref{sec4:eq:Cntrl_BacksteppingTrafo1} with respect to space and time and inserting the obtained derivatives,~\eqref{sec4:eq:Cntrl_BacksteppingTrafo2} and~\eqref{sec3:eq:ControlDesignModel_w123PDE} in~\eqref{sec4:eq:Cntrl_TargetSystem_a123PDE} yields the expressions $C^-(x,\xi)$ and $C^+(x,\xi)$. Finally, inserting~\eqref{sec3:eq:ControlDesignModel_bcL} and~\eqref{sec4:eq:Cntrl_TargetSystem_BCL} in~\eqref{sec4:eq:Cntrl_BacksteppingTrafo2} evaluated at $x=L$ and formulating the result with respect to the original physical variables, i.e. the densities and velocities of both vehicle classes yields
\begin{align}
&U(t) = -\kappa e^{\frac{\hat{J}_{44}}{\lambda_4}L}\bar{R}_1T_u^{-1}(L)\Psi(L,t) \notag \\
&-\kappa e^{\frac{\hat{J}_{44}}{\lambda_4}L}\hspace{-3pt}\int_0^L\hspace{-5pt}\left(K(L,\xi)T_u^{-1}(\xi)+L_{11}(L,\xi)T^{-1}_l(\xi)\right)\Psi(\xi,t)d\xi \label{sec4:eq:ControlLaw_PhysicalVariables}
\end{align}
with
\begin{align}
\Psi(\Omega,t) = z(\Omega,t)-z^*, \quad \Omega\in\{\xi,L\},
\end{align}
determining the control input. Notice that the transformation~\eqref{sec3:eq:SummarizedTransformation} is separated in two parts
\begin{equation}
\label{sec3:eq:Seperation_Tinv}
\left[\begin{array}{c}
T_u^{-1}(x) \\
T_l^{-1}(x) \\
\end{array}\right]= \left[\begin{array}{cccc}
0 & e^{-\frac{\hat{J}_{22}}{v_2^*}x} & 0 & 0 \\
0 & 0 & e^{-\frac{\hat{J}_{33}}{\lambda_3}x} & 0 \\
e^{-\frac{\hat{J}_{11}}{v_1^*}x} & 0 & 0 & 0 \\
0 & 0 & 0 & e^{-\frac{\hat{J}_{44}}{\lambda_4}x}
\end{array}\right]
\Theta^{-1},
\end{equation}
where $T_u^{-1}(x)\in \mathbb{R}^{3 \times 4}$ and $T_l^{-1}(x) \in \mathbb{R}^{1 \times 4}$. \\
In the following, the presented results are summarized in a theorem.
\begin{thm}
Consider the linearized two-class AR model~\eqref{sec2:eq:LinearizedModelEq} and~\eqref{sec2:eq:BC_lin}. Notice that the first vehicle class is assumed to represent faster vehicles with equilibrium speeds being higher than the ones of the second vehicle class. Assuming that the investigated traffic is in congested regime and the initial traffic density and velocity profiles satisfy
\begin{equation}
\tilde{\rho}_1(x,0),\tilde{v}_1(x,0), \tilde{\rho}_2(x,0),\tilde{v}_2(x,0) \in \mathcal{L}^\infty([0,L]),
\end{equation}
applying the control law~\eqref{sec4:eq:ControlLaw_PhysicalVariables} in~\eqref{sec2:eq:BCL_lin_Flow_wContrl} yields convergence of the density and velocity perturbations to the equilibrium at zero
\begin{equation}
\tilde{\rho}_{e,1}(x)\equiv\tilde{v}_{e,1}(x)\equiv\tilde{\rho}_{e,2}(x)\equiv\tilde{v}_{e,2}(x)\equiv 0
\end{equation}
in finite time $t_F$ given by~\eqref{sec4:eq:FiniteConvergenceTime_tF}. The kernel $K(x,\xi)$ is the solution of the well-posed kernel equations~\eqref{sec4:eq:kerneleq} and $L_{11}(x,\xi)$ follows according to~\eqref{sec4:eq:kerneleq_PDEL_Linserted}.
\end{thm}
\vspace{-7pt}
\section{Numerical Simulation}
\begin{figure*}[htbp]
	\begin{center}
		\includegraphics[width=14cm,height=5cm, trim = {0cm 1cm 0cm 2cm},clip]{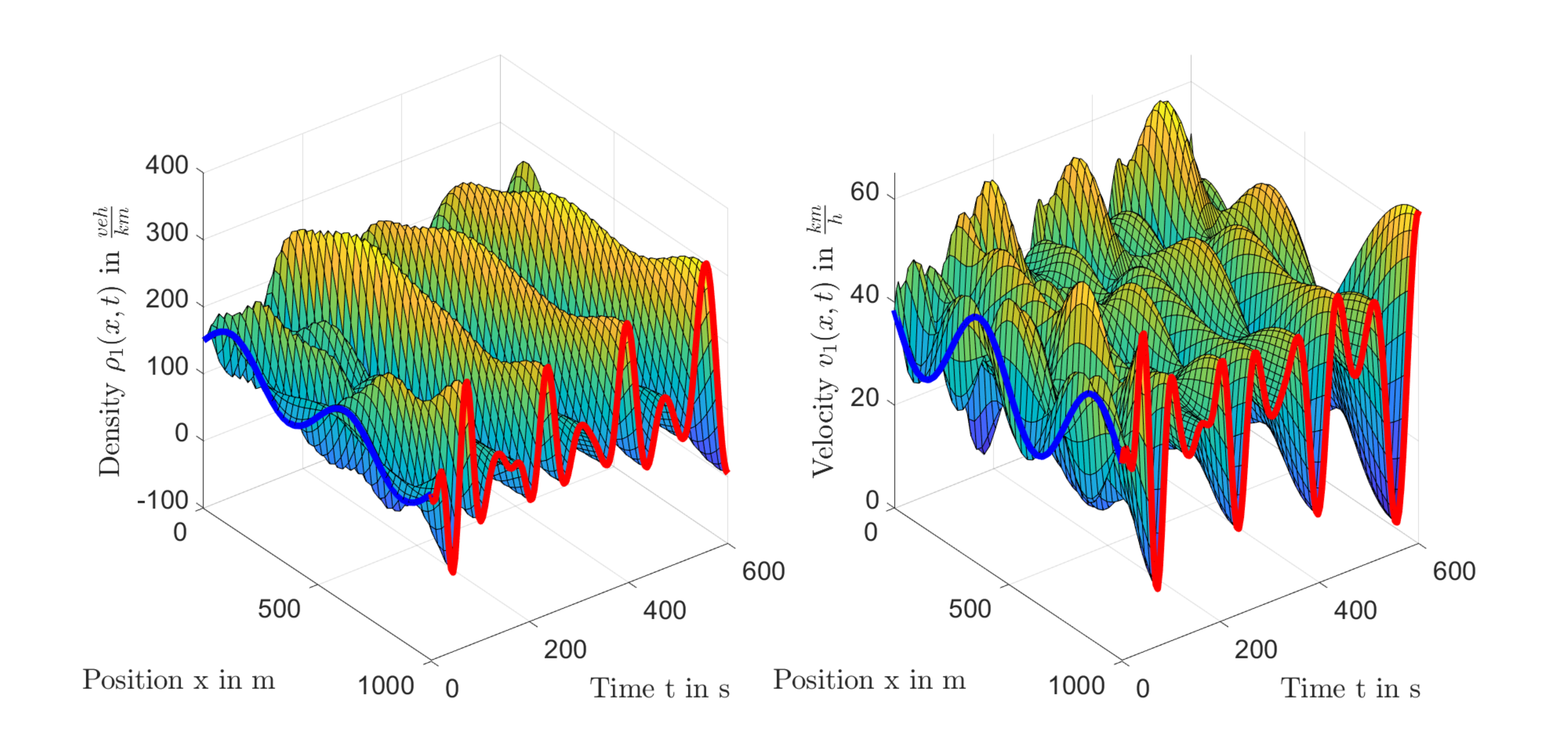} 
	\end{center}
\vspace{-7pt}
\caption{Traffic density and velocity of class $1$ without control.}
\label{sec6:fig:Class1_WithoutControl}
\vspace{0pt}
\end{figure*}

\begin{figure*}[htbp]
	\begin{center}
		\includegraphics[width=14cm,height=5cm, trim = {0cm 1cm 0cm 2cm},clip]{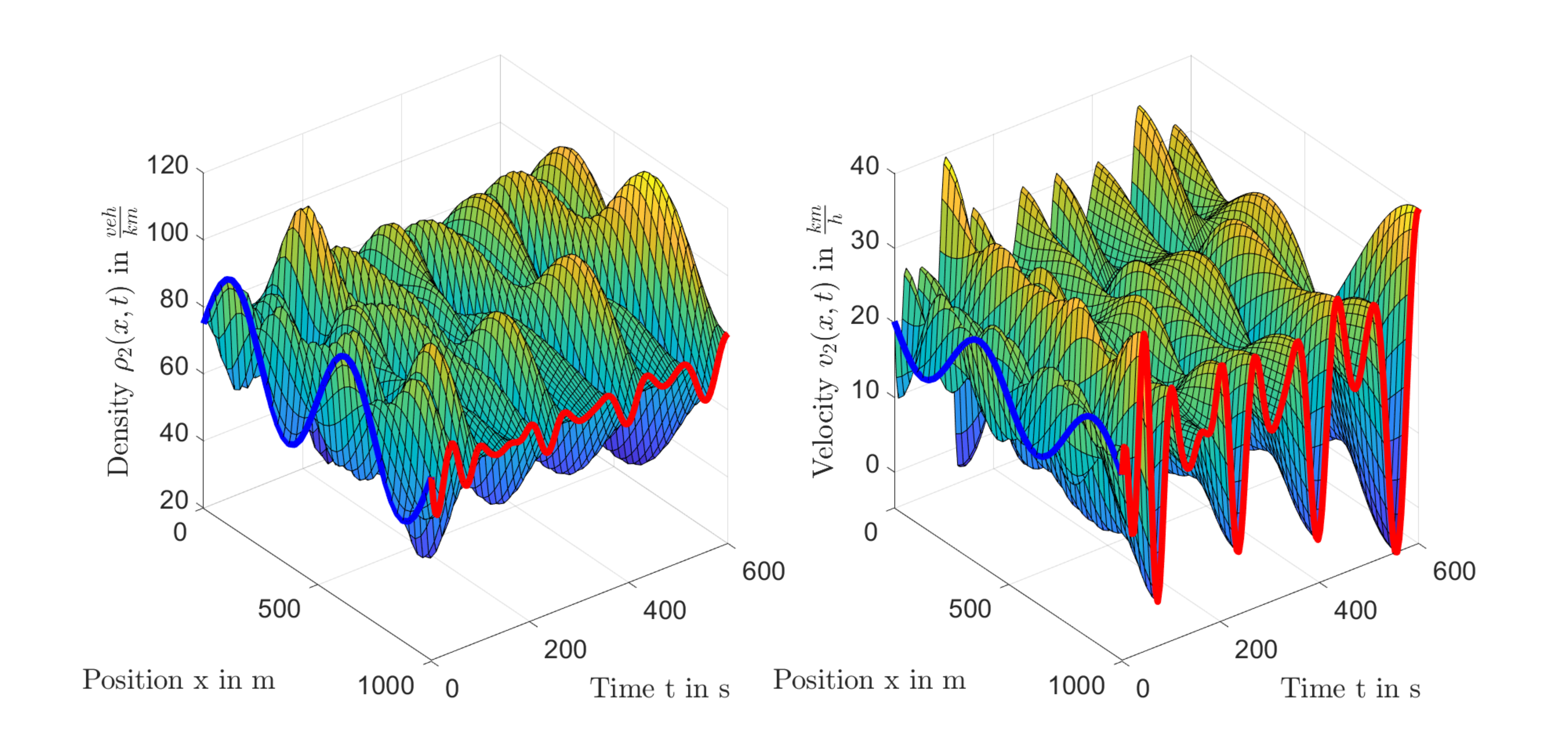} 
	\end{center}
\vspace{-7pt}
\caption{Traffic density and velocity of class $2$ without control.}
\label{sec6:fig:Class2_WithoutControl}
\vspace{0pt}
\end{figure*}

\vspace{-5pt}
The performance of the developed full state feedback controller is verified by a numerical simulation of the linearized model in the next step. All introduced parameters and the equilibrium densities $\rho_1^*$ and $\rho_2^*$ are chosen such that $\lambda_4<\lambda_2<\lambda_3<\lambda_1$ holds for the characteristic speeds. Hence, the traffic is in congested regime. Moreover, the initial profiles are given by
\begin{align}
\label{sec6:eq:Simulation_InitialProfiles}
\rho_i(x,0)\hspace{-2pt}=\hspace{-2pt}\rho_i^*\hspace{-2pt}+\hspace{-2pt}\frac{\rho_i^*}{4}\hspace{-2pt}\sin\hspace{-1pt}\left(\frac{4\pi}{L}x\right),\text{ }\hspace{-2pt}v_i(x,0) \hspace{-2pt}= \hspace{-2pt} v_i^*\hspace{-2pt}-\hspace{-2pt}\frac{v_i^*}{4}\hspace{-1pt}\sin\hspace{-2pt}\left(\frac{4\pi}{L}x\right)
\end{align}
and represent stop-and-go traffic since they describe alternating areas of dense, slow traffic and light, fast traffic. \\
Figure~\ref{sec6:fig:Class1_WithoutControl} shows the obtained open loop simulation results of the linearized two class AR model for density and velocity of vehicle class $1$, whereas Figure~\ref{sec6:fig:Class2_WithoutControl} depicts the results for vehicle class $2$. The initial profiles are highlighted with a blue line. The densities and velocities at the outlet, i.e. $\rho_i(L,t)$ and $v_i(L,t)$, are marked in red. The figures show that the initial stop-and-go oscillations are amplified. Furthermore, the closed loop simulation results are presented in Figure~\ref{sec6:fig:Class1_WithoutControl} and Figure~\ref{sec6:fig:Class2_WithoutControl}. The green line indicates the finite convergence time $t_F\approx \unit[237]{s}$. All states converge to their equilibrium values as suggested by the theory.
\begin{figure*}[htbp]
	\begin{center}
		\includegraphics[width=14cm,height=5cm, trim = {0cm 1cm 0cm 2cm},clip]{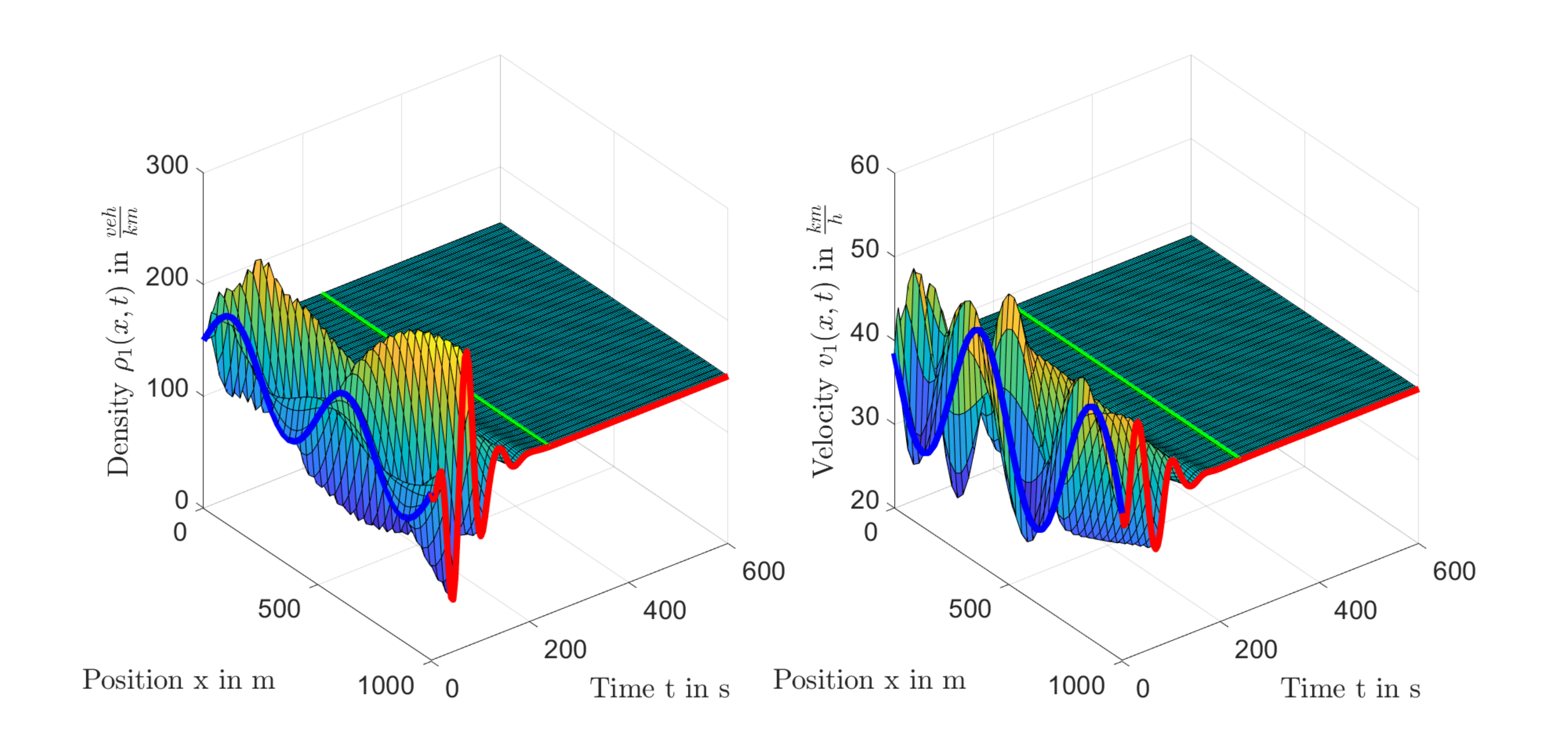} 
	\end{center}
\vspace{-7pt}
\caption{Traffic density and velocity of class $1$ with full state feedback control. The green line indicates $t_F$.}
\label{sec6:fig:Class1_Control}  
\vspace{0pt}
\end{figure*}
\begin{figure*}[htbp]
	\begin{center}
		\includegraphics[width=14cm,height=5cm, trim = {0cm 1cm 0cm 2cm},clip]{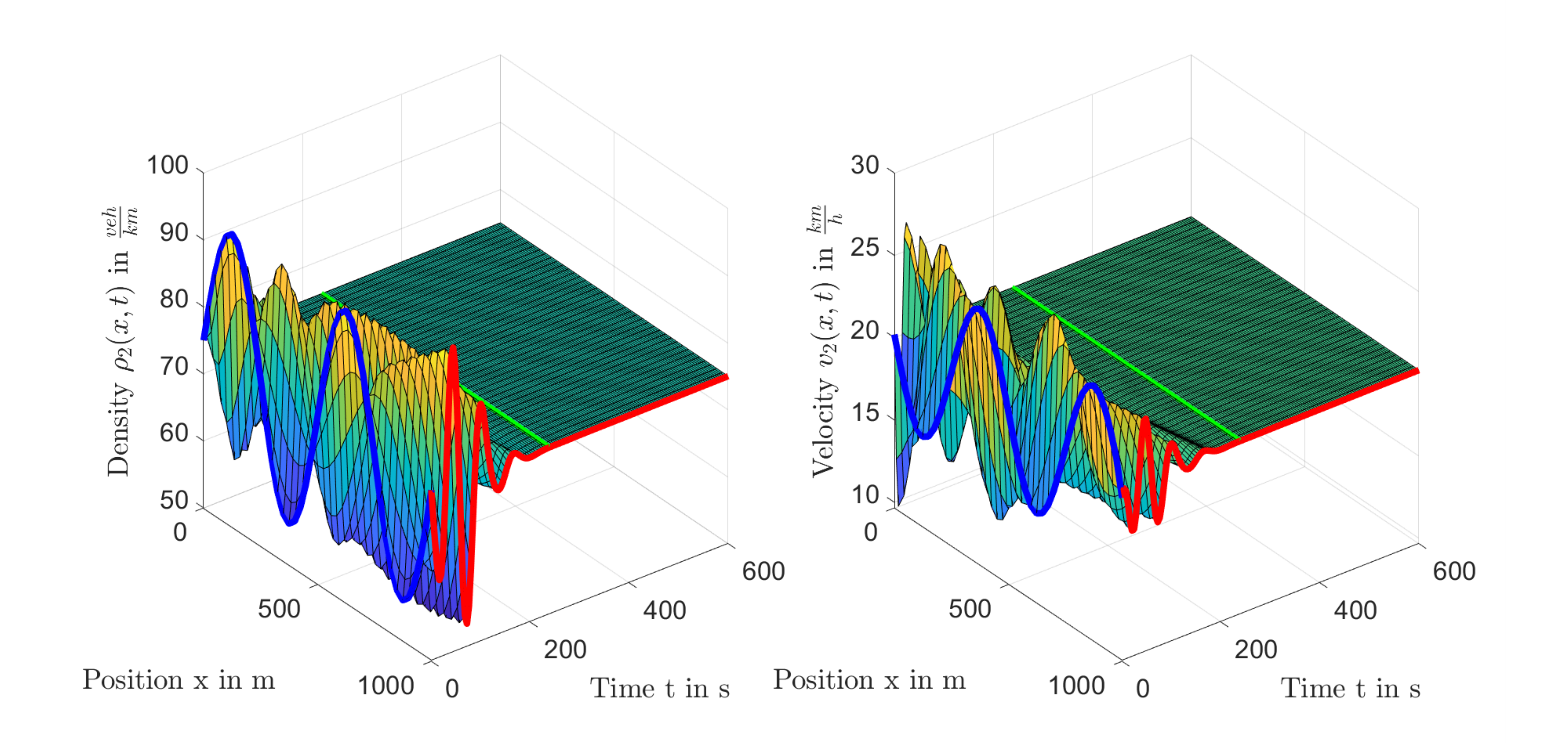}
	\end{center}
\vspace{-7pt}
\caption{Traffic density and velocity of class $2$ with full state feedback control. The green line indicates $t_F$.}
\label{sec6:fig:Class2_Control}  
\vspace{0pt}  
\end{figure*}
\vspace{-5pt}
\section{Concluding remarks}
\vspace{-5pt}
The presented feedback control result is only a first step of control design for multi-class traffic models using backstepping. A first future research topic is given by the extension of the results to traffic that is subdivided in three classes, e.g. distinguishing between motorcycles, average vehicles and trucks. The corresponding three-class AR model then consists of six nonlinear first order hyperbolic PDEs. Consequently, the signs of the characteristic speeds in the congested regime and their relation to the two-class result need to be investigated for that case. \\
Typically, it is hard and expensive to measure the entire state vector at every point along the investigated track section leading to a second future research topic: the full state feedback controller can be extended to an output feedback controller. Hence, the state vector is only measured at a single point and an observer generates estimates of the state vector at the remaining points. Afterwards, the developed control law of this work is reformulated based on these estimates.
\vspace{-2pt}
\begin{ack}
\vspace{-5pt}
Mark Burkhardt would like to thank Professor Oliver Sawodny for organizing the collaboration with Huan Yu and Miroslav Krstic. Furthermore, he acknowledges Kevin Schmidt for his help within fruitful discussions while this work was carried out. 
\end{ack}
\vspace{-2pt}
\bibliography{references}          
\end{document}